\documentclass[12pt]{article}
\usepackage{amsfonts}
\usepackage{amssymb}
\usepackage[usenames]{color}
\usepackage{epsfig}

\begin{document}

\input epsf.sty

\centerline{\bf S.P.Novikov\footnote{Sergey P. Novikov, IPST and
MATH Department, University of Maryland, College Park MD, USA and
Russian Academy of Sciences (Steklov Math Institute and Landau
Institute for Theoretical Physics, Moscow), e-mails
novikov@ipst.umd.edu, snovikov@mi.ras.ru; This work is partially
supported by the Russian Grant in the Nonlinear Dynamics. }}

 \vspace{0.3cm}

\centerline{\Large  New Discretization of Complex Analysis:}

 \centerline{\Large The Euclidean and Hyperbolic Planes}

\vspace{0.5cm}

{\it Abstract. Discretization of Complex Analysis on the Plane based
on the standard square lattice was started in 1940s. It was
developed by many people and also extended to the surfaces
subdivided  by the squares. In our opinion, this standard
discretization does not preserve well-known remarkable features of
the Completely Integrable System. These features certainly
 characterize the standard
Cauchy Continuous Complex Analysis. They played a key role in the
great success of
 Complex Analysis in Mathematics and Applications. Few years ago we
developed jointly with I.Dynnikov a New Discretization of Complex
Analysis (DCA) based on the two-dimensional manifolds with colored
black/white triangulation (see \cite{ND}). Especially deep results
were obtained for the Euclidean plane with equilateral triangle
lattice in the works \cite{ND,GN}. Our approach preserves a lot of
features
 of Completely Integrable Systems.
In the present work we develop  a DCA theory for the analogs of
equilateral triangle lattice in  the Hyperbolic plane. This case is
much more difficult than Euclidean. Many problems (easily solved for
the Euclidean Plane) are not solved here yet. Some specific very
interesting ''dynamical phenomena'' appear in this case: for
example, description of boundaries of the most fundamental geometric
objects (like the round ball) leads to dynamical problems. Mike
Boyle from the University of Maryland helped me to use here the
methods of symbolic dynamics. }

\vspace{2cm}

{\bf Introduction. History.} We do not discuss here ''geometric''
discretizations of conformal mappings started in early XX
century-see survey and history in \cite{Mer}. By the way, specific
topological properties of surfaces subdivided by the squares, was
started by the combinatorial geometry/topology group in the Steklov
Institute in 1980s by the suggestion of the present author (see
\cite{DSS}). This group called this object ''The Quadrillage''. I
formulated these problems under the influence of the Statistical
Physics of Lattice Models, after discussions with A.Polyakov. It is
mentioned also in my Encyclopedia volume number 13, Topology-I.--see
\cite{N2}. In the works of the group \cite{Mer} made in last decade
such surfaces subdivided into squares are called ''quad-graphs''.

{\bf Our goal is to discretize Cauchy-Riemann operator
$\bar{\partial}$ as a Linear Difference First Order Operator on
Triangulated Surfaces.} Let us say at the very beginning of this
 article that {\bf Neither in the Standard Approach nor in our New Discretization
 The Discrete Analogs of Holomorphic Functions form a Commutative Associative Ring
 (at least naturally).}

 The Discretization of Complex Analysis as a Linear Difference Operator
  was done for the square lattice in $R^2$ by
Lelong-Ferrand. Her work was published in 1944 (see \cite{LF}).  In
the work \cite{DUF} Duffin developed or rediscovered this approach.
He also extended it to the rombic lattices later.

{\it Duffin claimed that
another mathematician already discovered this idea before and published his work
in 1941 in the notes of some small provincial South American University.
 Our attempts to
find this obscure edition failed, so we continue to quote Lelong-Ferrand as a
first author
 who invented this idea.}

 In the standard approach  Discrete Analog  of the Cauchy-Riemann Operator
$\bar{\partial}$ acts on the $C$-valued functions $\psi$ of vertices
in the square lattice on the plane:
$$Q_{square}\psi(m,n)=\psi(m,n)+i\psi(m+1,n)-\psi(m+1,n+1)-i\psi(m,n+1)$$
By definition, every solution to the equation $Q_{square}\psi=0$ is
a d-holomorphic function.

A lot of people work here now developing this approach (see, for
example,  in \cite{Mer}).

{\it Let us present  a hint of the idea: Why we are not satisfied by
the standard
 discretization? Of course, every continuous model admits many  difference
  approximations which converge to it in the continuum limit.

  For the continuous systems with huge hidden algebraic symmetry
  we wish to find optimal discretization preserving as much  symmetry as possible.
  As we already know, there exists no natural discretization preserving multiplication
  of holomorphic functions. All attempts to invent multiplicative structure with
  good properties failed. It became clear many years ago.}

Let us point out following:

1. In the standard approach discrete analog of the Cauchy-Riemann
operator is in fact a
 second order difference operator:  two lengths  are
involved in the sum (the lengths of  sides and  diagonals).

 {\bf We wish to have
 a first
 order difference operator as a natural discrete analog of $\bar{\partial}$.
 Is it possible? The answer is YES for the Equilateral Triangle Lattice. Corresponding theory
  of the difference first order ''Triangle Operators'' on the simplicial complexes
  is developed in the Part 1 in more general form. We define also nonstandard
  discrete analogs of
  $GL_n$-Connections for the triangulated n-manifolds}.

  In the Part 2 we develop a general theory of discrete analogs of holomorphic
   functions
(the d-holomorphic functions) on the 2-manifolds with discrete
analog of conformal structure (DCS). {\bf We define DCS as a colored
black/white triangulation.}

The d-analogs of Liouville Principle and Maximaum Principle will be clarified here
 for the general 2-surfaces with DCS.

 Part 3 is dedicated to the theory of d-holomorphic functions on the Equilateral
 Triangle Lattice.

2.Continuous 2D Laplace operator admits a natural factorization
$$\Delta=\bar{\partial}\partial$$ Nothing like that exists on the square lattice.

{\it However, Equilateral Triangle Lattice is much better for this--see below.}

{\bf Problems:}

1.How to construct Discrete Analogs of Holomorphic Polynomials without multiplication?

People did it with great efforts after many years in the standard approach. Their constructions
are non-
canonical and non-unique.

{\bf It is very easy in our new approach.
Our construction is obvious and canonical.}

2.How to construct rational functions without multiplication?

 In particular, in the standard continuous complex analysis rational functions
 are especially important because the function $1/z$ (the Cauchy Kernel) is also a
unique fundamental solution to the equation $$\bar{\partial}(\psi(z))=2\pi i\delta(z)$$
decreasing at infinity?

{\bf The answer to this question is also positive and canonical in our approach.
It was missed in our first
 work with Dynnikov \cite{ND} and made by Grinevich and R.Novikov later in \cite{GN}.}

Part 4 is dedicated to the Hyperbolic (Lobachevski) Plane $H$.

 Equilateral Triangle Lattice in $H$ is defined simply as a triangulation of plane such that
  $p$ triangles hit every vertex and
$p\geq 7$. We need $p=2s$. For all $s\geq 4$ we have a triangulation of Lobachevski Plane $H$
such that it admits
a black/white coloring. Let us concentrate our attention on the case $s=4$.

 At the moment we can
construct neither d-analog of polynomials nor d-rational functions.

 Even boundary of
round ball $D_r$ of radius $r\in Z$ is a complicated object. How to
describe it? How many points does it contain? We succeeded in the
solution of this problem using methods of symbolic dynamics.
 Mike Boyle
from the University of Maryland helped us. The boundary $\partial
D_r$ contains approximately $\lambda^r$ points where
$\lambda=2+\sqrt{3}$.

It turns out that for the reconstruction of d-holomorphic function
in the ball $D_{r+1}$ we need to know  values of the boundary
function in $[n(r+1)/2]+1$ points. Here $n(r)$ is a number of
vertices in the boundary. It exactly corresponds to the continuous
case where only set of coefficients at the nonnegative powers of $z$
(i.e. one half plus one Fourier coefficients of function on the
boundary circle) determines completely our function inside. The main
unsolved problem is:

{\bf Problem. How to find a bounded basis of d-holomorphic functions
on Hyperbolic Plane?}

\vspace{1cm}

{\bf Part I. Definitions. Discrete $GL_n$-Connections. B/W
manifolds. The First Order Triangle Operators.}

\vspace{0.5cm}

Take a simplicial complex $M$. We assume that this complex is given
with canonical metric where every simplex is a standard unit linear
subsimplex in euclidean space with natural euclidean metric. We fix
a family of $n$-simplices such that every vertex belongs at least to
one simplex of that family. Fix also collection of nonzero numbers
$b_{T:P}\neq 0$  for all $T\in X$ and $P\in T$.

\newtheorem{defi}{Definition}
\begin{defi} We call following operator
 the Triangle Operator $Q^X$ associated with
family $X$: $$Q^X\psi(T)=\sum_{P\in T}\psi (P)b_{T:P}$$ It maps
functions of vertices into functions of simplices $T\in X$. Our
coefficients normally belong to $R$ or $C$. The operators $Q^X$ are
the first order linear difference operators.
\end{defi}

Now we define a Discrete Version of $GL_n$-Connections.

\begin{defi}
Let $X$ be a family of all $n$-simplices in $n$-manifold $M$. We
call the equation $Q\psi=0$ Discrete Differentially-Geometrical
$GL_n$-Connection. The coefficients $b_{T:P}$ are defined in every
$n$-simplex up to nonzero factor. So every DG-Connection is defined
by the set of ratios $$\mu^T_{PP'}=b_{T:P}/b_{T:P'}$$
\end{defi}

 This
discretization is different from the standard Wilson Discretization
used by physicists studying Yang-Mills fields; There is no natural
way to select compact holonomy groups in our approach.

 The theory of
discrete $GL_n$ connections was constructed in the works
\cite{DN,N}, some first ideas appeared in \cite{N1,ND} in connection
with completely integrable systems. We are working only with
scalar-valued functions. At the same time we are going to define a
Nonabelian Curvature and Holonomy with values in the group $GL_n$
for the $n$-manifolds.

As we  demonstrated in the \cite{N}, two different Holonomy
Representations can be constructed here: Abelian (''Framed
Abelian'') and Nonabelian. We need here only Nonabelian Holonomy.

What is  a Nonabelian Discrete Curvature? A Nonabelian Holonomy
group is defined along the Thick Paths.

\begin{defi}
 We call by Thick Path any sequence of $n$-simplices $<T_1T_2...T_k>$
 such that intersection $T_j\bigcap T_{j+1}$ is exactly $n-1$-dimensional
 face $\Delta_j$, and $\Delta_j\neq\Delta_{j+1}$. We say that  Thick Path is closed if
 $T_k= T_1$. Its length is equal to $k-1$.
\end{defi}

Therefore in every simplex $T_j$ of the thick path exactly two
different  $n-1$-dimensional faces are fixed: $\Delta_j\subset
T_{j+1}$ is called an ''In-face'' in $T_{j+1}$, and $\Delta_J\subset
T_j$ is called an ''Out-face'' in $T_j$.

Let us define a {\bf Parallel Transport} of the $n$-vector $\psi$
consisting of values in all vertices of the in-face $\Delta_1\subset
T_1$, along the thick path $<T_1...T_k$. Solving equation $Q=0$, we
define this function in all vertices of $T_1$. So we know values of
$\psi$ in all vertices of the out-face in $T_1$ equal to the in-face
in $T_2$. After that we are doing the same procedure for $T_2$ and
so on. Finally we define a ''Transported value'' of $n$-vector
$\psi$ in the out-face of the simplex $T_k$. This map is linear. For
the closed paths we are coming to the {\bf Nonabelian Holonomy
Representation of the Semigroup of the closed thick paths into the
Group $GL_n$.}

$$\Omega_{thick}(M,\Delta):\rightarrow GL_n$$

\begin{defi} We call by the Nonabelian Discrete Curvature in the vertex $P$  a Holonomy
Linear Map along the Thick Closed Paths in the Simplicial Star
$St(P)$. It is enough to know holonomy for all thick closed paths
around all $n-2$=simplices--see Fig 1 for $n=2$.

\end{defi}

\begin{center}
\color{black}Fig 1 \vspace{5mm}
 \mbox{\epsfysize=4cm\epsffile{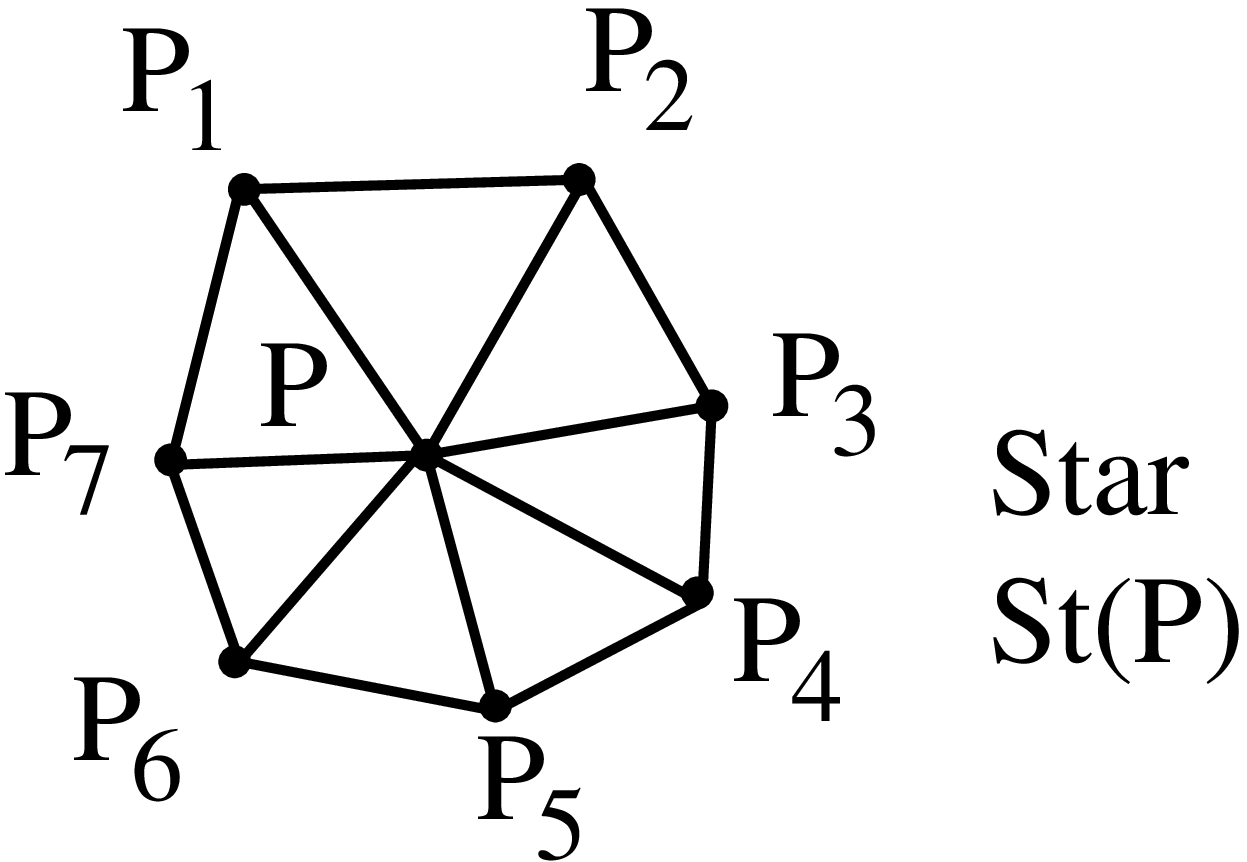}}
\end{center}

However, the natural gauge group
$$\psi\rightarrow f\psi, Q\rightarrow g(T)Qf^{-1}(P)$$ is abelian,
$f\neq 0,g\neq 0$. There are very interesting specific features here
realizing non-abelian $GL_n$ Connections using only the spaces of
scalar functions of vertices. Some sort of mixing abelian and
non-abelian properties appears. Details of classification can be
found in the work \cite{N}. No theory of Characteristic Classes is
constructed yet.

\begin{defi}
We call special discrete $GL_n$-Connection with $b_{T:P}=1$, a
Canonical Connection
\end{defi}

     Consider now any $n$-manifold. Is it possible to color its $n$-simplices
      by the
      black and white colors
such that for every pair of $n$-simplices attached to each other
along $n-1$-face, the colors are opposite?

\begin{defi} We call such colored manifolds a B/W
manifolds.
\end{defi}

Such coloring exists if and only if every closed thick path consists
of even number of $n$-simplices. More general, fix any $n-1$-simplex
$\Delta$. We have a natural homomorphism
$$\Phi_2:\Omega_{thick}(M,\Delta)\rightarrow Z_2=Z/2Z$$
defined as a number of $n$-simplices in any thick path modulo 2.
This map  also looks like  simplest Connection whose Holonomy Group
is $Z_2$, but we do not realize it by the linear operator.

\begin{defi} Let the map $\Phi_2$ is trivial for
 all thick closed paths in the star of every vertex. We call such triangulation
 ''Locally B/W''.
\end{defi}

   In this case following
  simple lemmas are true:

\newtheorem{ex}{Example}

\newtheorem{lem}{Lemma}

\begin{lem}
Under this assumptions the map $\Phi_2$ defines a map $\phi_2:\pi_1(M)\rightarrow Z_2$.
\end{lem}

\begin{lem}
Curvature of the Canonical Connection is trivial for every locally B/W
manifold.
\end{lem}

Proof. Fix any $n-2$-simplex $S$ in the star $St(P)$. Let exactly
$s(S)$ $n$-simplices contain $S$. The Nonabelian Holonomy with
initial values $a_0,...,a_n$ in the vertices of $n$-simplex
$T_n=<S,P_{n-1},P_n>$ containing $S$ (and around it), is determined
by permutation $A^{s(P)}$. Here $A(a_0,...,a_{n-2},a_{n-1},a_n)=
(a_0,...,a_{n-2},a_n,a_{n-1})$. So $A^2=1$. We see that  Nonabelian
Holonomy (Curvature) is equal  to $1$ if and only if $s(S)$ is even.
Lemma is proved.

Every locally flat Connection (i.e. Nonabelian Curvature is trivial
for every vertex) defines a homomorphism of topological holonomy
$\pi_1(M)\rightarrow GL_n$.

\begin{lem}
 Every flat  Canonical Connection  defines a Topological
Holonomy Homomorphism. Its image  belongs to the subgroup
$S_{n+1}\subset GL_n$.

\end{lem}

Proof. The imbedding is following: Take arbitrary distinct $n+1$
numbers $a_0,...,a_n$ such that $\sum a_j=0$. Put this set of
numbers in the vertices of initial $n$-simplex $T_1$. For every
thick closed path $T_1,....,T_k$ the values of Parallel Transport
along this path in every simplex $T_j$ consists of the very same
numbers $a_q$ because $\sum a_j=0$, and Connection is Canonical. In
the final simplex $T_k=T_1$ we obtain permutation belonging to
$S_{n+1}$ which does not depend on the choice of initial numbers
$a_j$. Lemma is proved.

{\bf Now we are going to consider only Orientable Manifolds $M$.}

\begin{lem} Let Curvature and Topological Holonomy are trivial for the
Canonical Connection (i.e. it is flat). Then B/W coloring exists
globally.
\end{lem}

Proof. Let a closed thick path $<T_1,...,T_k>$ defines a trivial
permutation corresponding to the Canonical Connection. We claim that
this path consists of the even number of $n$-simplices. Indeed, it
follows from the fact that our Nonabelian Holonomy is equal to the
product of elementary transpositions. Every one of them changes
orientation. The whole path preserves orientation. So our lemma
follows.

\begin{lem}
For every manifold $M$ with flat Canonical Connection the equation
$Q\psi=0$ defines  $n$-dimensional linear space of Covariant
Constants.
\end{lem}

Proof obviously follows from the fact that every initial values in
the $n-1$-simplex $Delta$ defines a covariant constant correctly
because our connection is flat.

So we fix a manifold $M$ with B/W structure. We always assume that
our Canonical Connection is flat. It is always true for simply
connected manifolds. It is locally flat because B/W structure is
given.For non-simply-connected case B/W structure implies that
Nonabelian Holonomy maps $\pi_1(M)$ into $S_{n+1}$. So Canonical
Connection became globally flat on some finite covering with number
of sheets no more than $n!$. We define d-analog of holomorphic
functions in every B/W manifolds (see below).

 {\bf How to construct Continuous Limit
for our Discretization of Complex Analysis?}

For that we are doing the following procedure: Let $n=2$. Construct
a special Covariant Constant $f_0$ whose values in every triangle
are $1,\zeta,\zeta^2$ where $\zeta^3=1$. It is unique up to the
group $S_3$ permuting vertices of the initial triangle. Apply
following gauge transformation to the Canonical Connection:
$$Q\rightarrow  f_0^{-1}Qf_0, \psi\rightarrow f^{-1}_0\psi$$
We still have two-dimensional space of covariant constants for this
connection but in this gauge  one of covariant constants became an
ordinary constant. We work over the extended field containing
$\zeta$. For the case $R$ and extended field $C$ we used this
construction to get  Continuum Limit of our theory to the ordinary
complex analysis. For $M=R^2$ exactly one half of our theory
converges to the ordinary complex analysis. The second half diverges
for small triangulations.

{\bf Now we start to construct New Discretization of Complex
Analysis.}

Let $M$ be any B/W-manifold. Following two families play fundamental
role in our theory: $X_1=b$ (all black simplices) and $X_2=w$ (all
white simplices)--see Fig 2.

\begin{center}
\color{black} Fig 2 \vspace{5mm} \mbox{\epsfysize=7cm
\epsffile{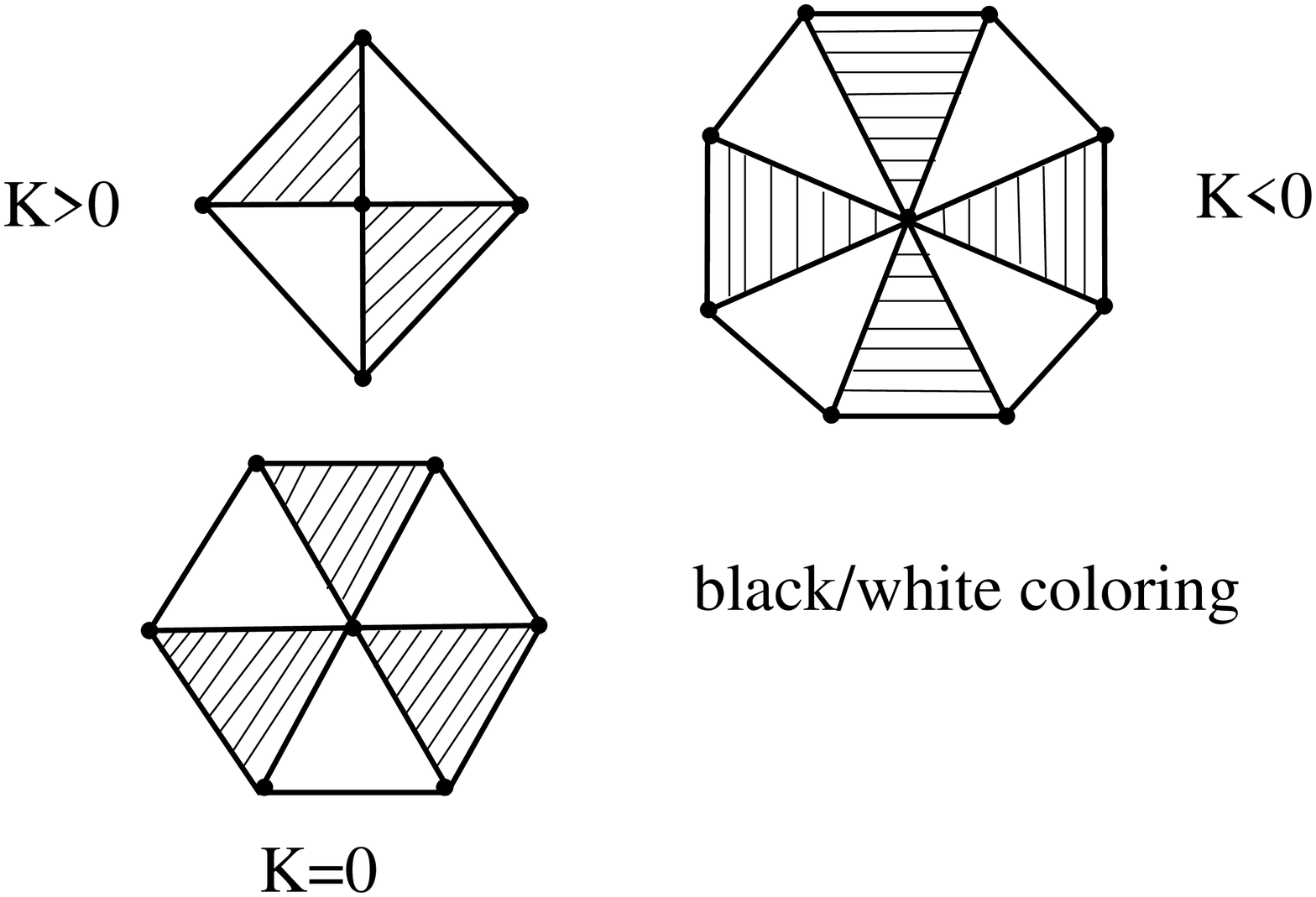}}
\end{center}

\begin{defi}The  operators $Q^b,Q^w$ are called Black and White
Triangle Operators. The union of this families $b\bigcup w$ is
complete. Using  the Canonical Connection, we define an operator
$Q=Q^b\bigoplus Q^w$.
\end{defi}

Following lemma is very simple:

\begin{lem} The operators $Q^b,Q^w,Q$ have following properties:
$$L= Q^*Q=2Q^{b*}Q^b=2Q^{w*}Q^w$$
Let  $\Delta=dd^*$ be an ordinary scalar simplicial Laplace-Beltrami
Operator, $m_P$ is   the number of  edges entering vertex $P$. For
2-manifolds we have $Q^*Q=-2\Delta+3m_P$.
\end{lem}

Proof follows from thew elementary calculation.

\begin{defi}Consider  any function $\psi$ such that $Q^b\psi=0$. For
$n=2$ we call such functions {\bf $d$-holomorphic}. We call
solutions to the equation $Q^w\psi=0$ $d$-antiholomorphic.
\end{defi}

Now we assume that the Canonical Connection is flat. For every black
triangle $T$ we define a unique covariant constant $E\psi(T)$ whose
values in $T$ coincide with $\psi$.

\begin{defi} We call $E\psi$ an {\bf Evaluation of function $\psi$ in the
black triangle $T$}.
\end{defi}

{\bf As a conclusion, we can view every   d-holomorphic function
either as a real-valued function $\psi$ on the vertices or as a
$R^2$-valued function $E\psi(T)$ of the black triangle $T$.}

\vspace{2cm}

{\bf Part II. D-holomorphic functions on 2-manifolds. Liouville
Principle and Maximum Principle}

\vspace{1cm}

As above, we consider triangulated orientable B/W manifolds $M$ with
flat Canonical Connection. Our field here is $R$. Following general
properties of d-holomorphic functions  were found in \cite{DN}:

\newtheorem{thm}{Theorem}
\begin{thm}(The Liouville Principle).
For the compact closed 2-manifold every d-holomorphic function is a
covariant constant.

\end{thm}

Proof is based on the ''instanton phenomenon'' for the quadratic
functional $(L\psi,\psi)$ if $L$ is factorizable $L=Q^*Q$. Let the
global minima $\psi$ in the Hilbert space $L_2$ is a zero mode
$L\psi=0$. The order of Euler-Lagrange equation drops twice reducing
to the ''Self-Duality Equation'' which is in our case simply the
d-analog of the Cauchy-Riemann Equation. Let us remind here that
{\bf This Property is the most fundamental feature of the standard
continuous Complex Analysis unifying it with The Completely
Integrable Systems}. The proof of our theorem is following:
$Q^b\psi=0$ implies $(Q^b)^*Q^b\psi=0 $ implies $L\psi=Q^*Q\psi=0$
implies $(L\psi,\psi)=0$ implies $(Q\psi,Q\psi)=0$ implies
$Q\psi=0$, so our theorem is proved.

 In
principle, the famous Instanton Phenomena discovered by Polyakov,
Belavin, Schwarz and Tiupkin in 1974 in the Nonlinear Variational
Calculus (Yang-Mills Field Theory), is based on the similar
arguments.

Let $D$ be a bounded domain in $M$ consisting of black triangles.
\begin{defi} We call $T\in D$ a boundary triangle if some of its
vertices belong to black triangle not belonging to the domain $D$.
\end{defi}

\begin{thm} (The Maximum Principle)
For every d-holomorphic function $\psi$ in $D$ the set of covariant
constants $E\psi(T)\in R^2$ for all $T\in D$ is contained in the
convex hull of the image of boundary triangles.
\end{thm}

Proof can be found in \cite{DN}.

In the next chapters we consider specific results obtained for the
Euclidean  (Part III) and Lobachevski or Hyperbolic Planes (Part
IV).

\vspace{1cm}

{\bf Part III. Equilateral Triangle Lattice in $R^2$.}

\vspace{1cm}

The standard regular equilateral lattice in $R^2$ (see Fig 3) admits
two basic (unitary) shift operators $t_1,t_2$ acting on vertices
$$t_1(m,n)=(m+1,n), t_2(m,n)=(m,n+1)$$.

\begin{center}
\color{black} Fig 3  \vspace{5mm}
 \mbox{\epsfysize=5cm\epsffile{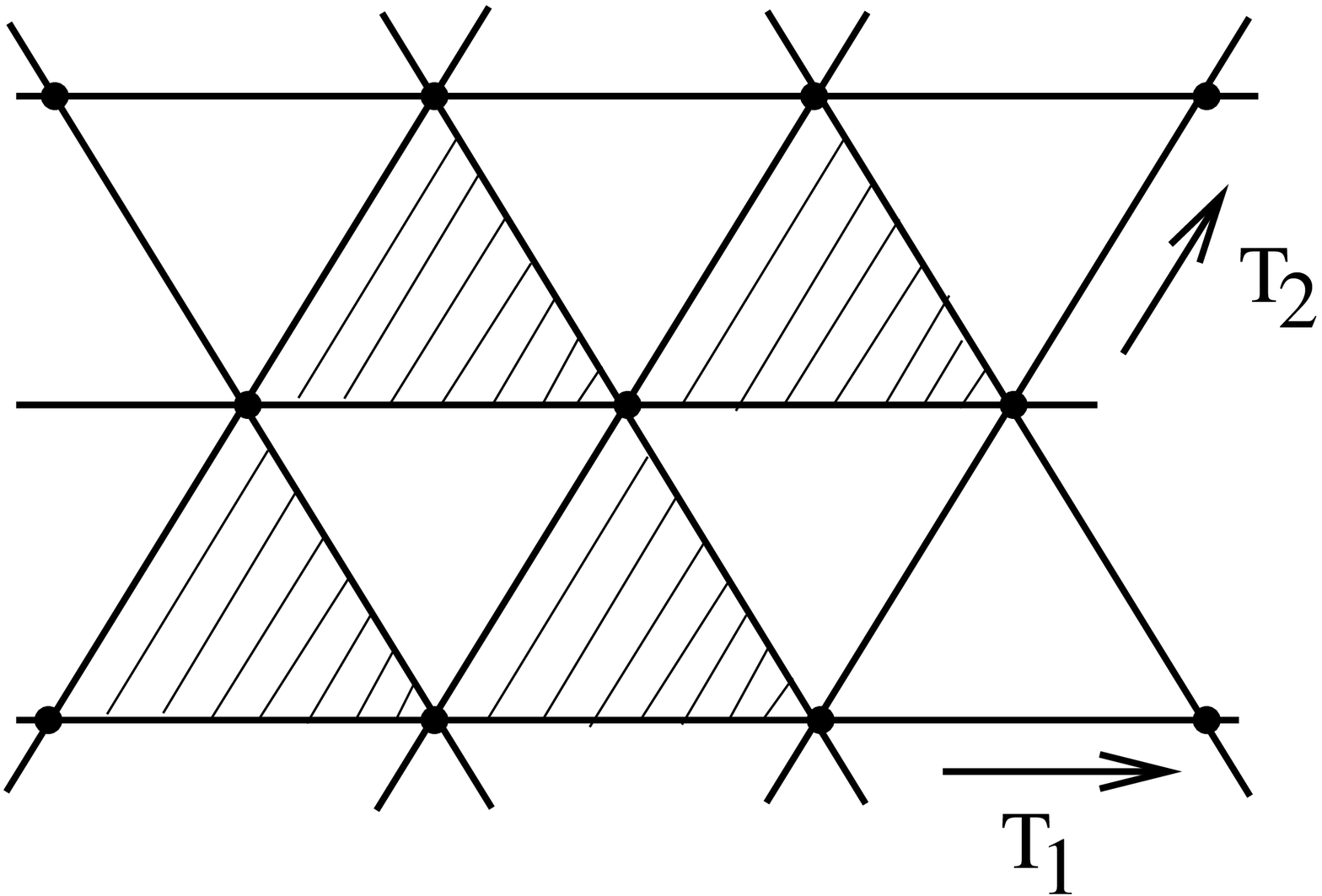}}
\end{center}

 All six elementary shifts
$t_i^{\pm 1}, (t_1t_2^{-1}))^{\pm 1}$ have equal length. The
operators $Q^b,Q^w$ are equal to
$$Q^b=1+t_1+t_2, Q^w=(Q^b)^*=1+t_1^{-1}+t_2^{-1}$$
They map the space of functions of vertices into itself and are
adjoint to each other. Our field  here is $R$. In the works
\cite{N1,ND} we defined first time  black and white triangle
operators and invented the idea of new type discrete $GL_2$
connection. We proved also that every difference second order
self-adjoint operator $L=a+bt_1+ct_2+dt_1t_2^{-1}+(adjoint)$
admits Laplace-type factorization $L=Q^*Q+W$ where $Q$ is some
triangle operator $Q=u+vt_1+wt_2$. and $W,a,b,c,d,u,v,w$ are some
real functions on this lattice. In particular, $-\Delta+9=Q^bQ^w$
for the standard Laplace-Beltrami operator. A theory of discrete
completely integrable systems based on the discretized  second
order operators, was started.

\begin{defi} We call d-holomorphic function $\psi$ polynomial (i.e. $\psi\in
Pol_k$)
of degree
$k$ if  $Q^b\psi=0$ and $(Q^w)^{k+1}\psi=0$.
\end{defi}

As it was established in  \cite{DN}, these functions have a
$k$-polynomial growth at infinity in $R^2$. They are completely
determined by their values in any  standard $k$-triangles $T_k$ ,
black from inside,  with $2k+2$ points in each edge (see Fig 4).

\begin{center}
\color{black} Fig 4 \vspace{5mm}
\mbox{\epsfysize=5cm\epsffile{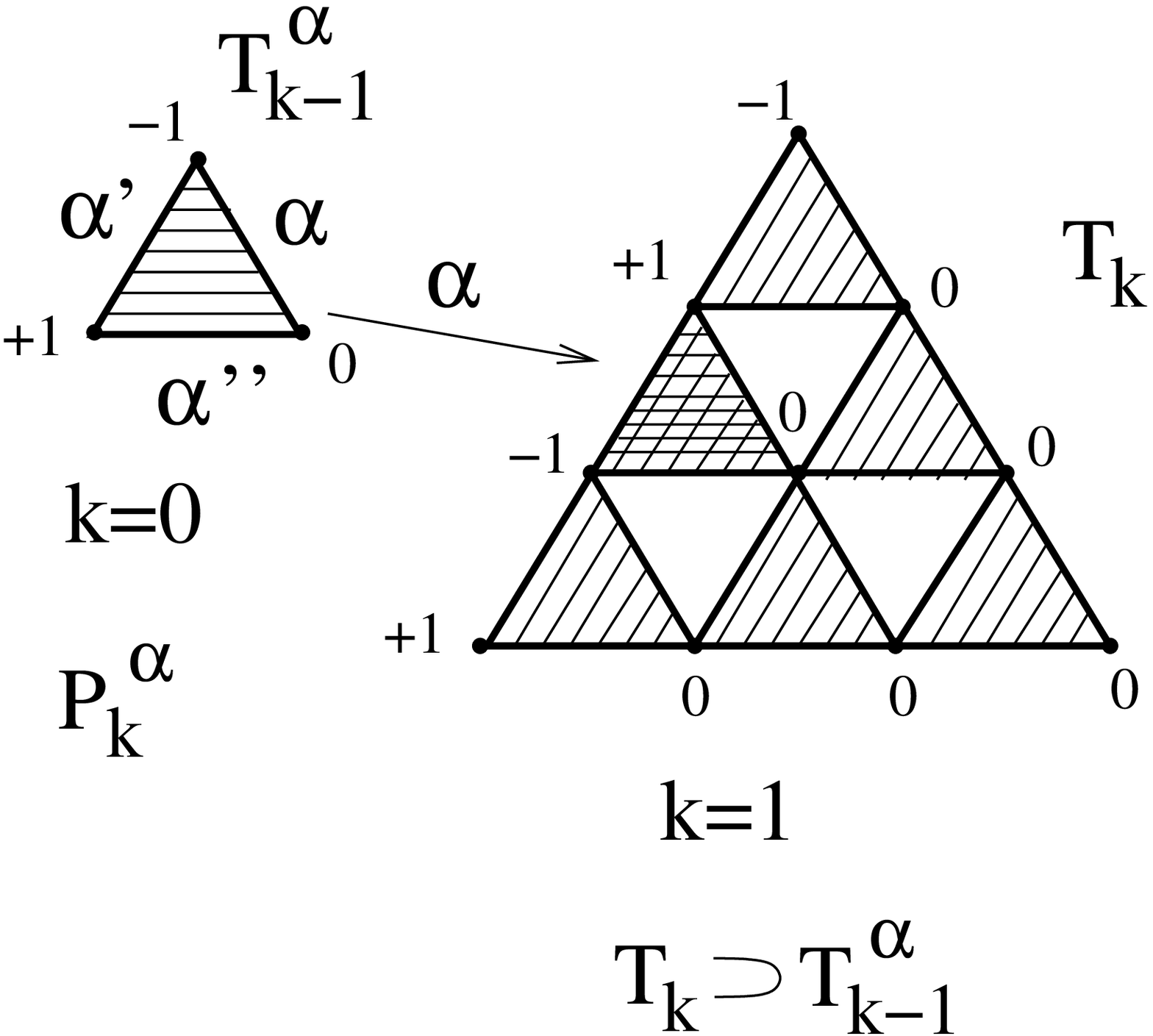}}
\end{center}

 So
the dimension of the space $Pol_k$ of $k$-polynomials is $2k+2$. We
can choose the canonical $k$-polynomials
$\psi_{T_k,\alpha},\alpha=1,2,3$. They are equal to zero in $T_k$
everywhere except one boundary edge $\alpha$ where they have values
$\pm 1$ . Sum of them belongs to the space $Pol_{k-1}$:
$$\sum_{\alpha}\psi_{k,\alpha}\in Pol_{k-1}$$
Therefore for the choice of basis we have to select pair of edges in
$T_k$.

 Following theorem (\cite{DN}) give natural analog of polynomial
approximation of holomorphic functions leading to  the ''d-Taylor
Series'':

\begin{thm} For every function $\psi$, d-holomorphic in $R^2$, and every
canonical triangle $T_k$, there exists a unique $k$-polynomial
$\phi$ such that $\psi-\phi$ is identically equal to zero in the
triangle $T_k$.

\end{thm}

For the Taylor decomposition of $\psi$ we have to choose an
increasing sequence of triangles $T_k,k=0,1,2,...$ and pair of
boundary edges in each of them. The choice of such sequence is
non=canonical.

{\bf How to get analog of the Cauchy formula?}

 We need to construct
a ''Cauchy Kernel'' (d-analog of $1/z$)  satisfying to equation:
$$Q^bG(x,y)=\delta(x-y)$$ where $x,y$ are points in the lattice
$x=(m,n),y=(m',n')$, and difference operator $Q^b$ acts on the
variables $x=(m,n)$. Having any such function, we consider any
bounded domain $D$ and $d$-holomorphic function $\psi$ in $D$. Let
us extend this function to the function $\bar{\psi}$ such that
$$\bar{\psi}(x)=\psi(x), x\in D; \bar{\psi}(x)=0, x\in R^2  minus D$$

We have $$\sum_y [Q^b\bar{\psi}(y)]G(x-y) =\psi(x)$$ for all $x\in
D$.  Let us point out that the function $Q^b\bar{\psi}(x)$ is
concentrated along the ''boundary strip'', so it is really analog of
the Cauchy formula. In the work \cite{DN} we constructed following
''hyperbolic-like'' Cauchy Kernel (see Fig 6). We call it Pascal
Triangle. It is equal to zero outside of infinite triangle and has
exponential growth inside of it.

\begin{center}
\color{black} Fig 5 \vspace{5mm}
 \mbox{\epsfysize=5cm\epsffile{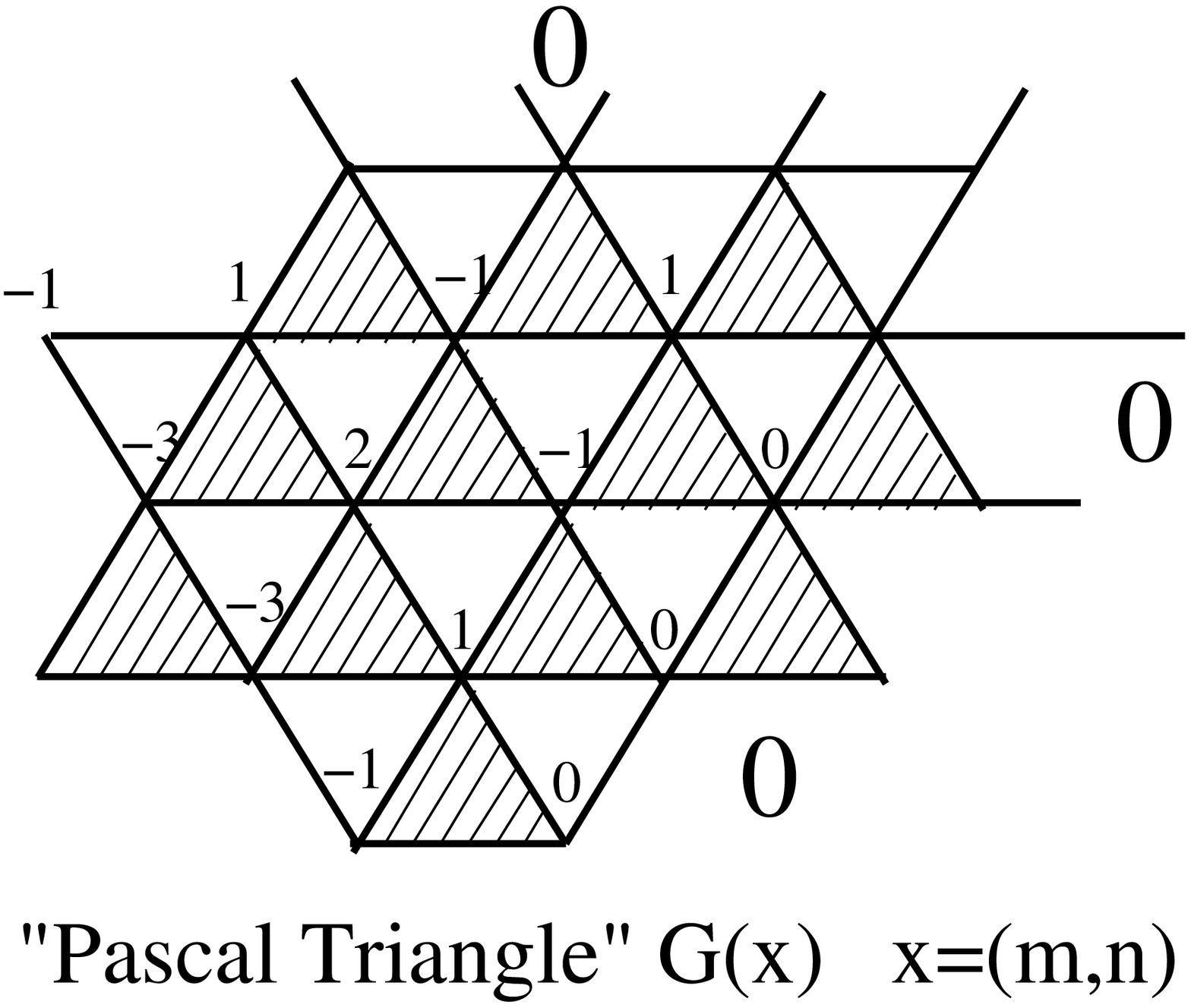}}
\end{center}

As it was pointed out by P.Grinevich and R.Novikov in \cite{GN},
much better Green function can be obtained simply by the Fourier
Transform
$$G(x)=(const)\int_0^{2\pi}\int_0^{2\pi}dk_1dk_2e^{imk_1}e^{ink_2}
(1+\exp\{ik_1\}+\exp\{ik_2\})^{-1}$$
where $x=(m,n)$.

They proved following theorem:

\begin{thm}This integral converges. This Green Function decays at infinity as
$c/d(P)$ where $d(P)$ is equal to the distance of this point to
zero. In particular, every d-holomorphic function on the whole
equilateral triangle lattice whose grouth at infinity is no more
that polynomial, is a d-holomorphic polynomial.

\end{thm}

{\bf Conclusion. This result gives a unique function decaying like
$1/|x|$ at $|x|\rightarrow\infty$. All functions $(Q^w)^kG(x)$ has a
growth like $1/|x|^{k+1}$. Therefore we have also analogs of
rational function without any  of multiplication (which does not
exist).}

\vspace{0.5cm}

{\bf Part IV. D-holomorphic functions on The Hyperbolic
(Lobachevski) Plane.}

\vspace{1cm}

Every triangulated plane $H^2$ such that more than six triangles
(edges) hit every vertex, can be viewed as a negative curvature
plane. An Equilateral Triangle Lattice we get in the case then this
number is the same for all vertices (i.e. $m_P\geq 7$). We need also
$B/W$ structure so our Equilateral Triangle Lattices $H^2_m$ in
$H^2$ with d-conformal structure (i.e. with B/W coloring) correspond
to $m=8,10,12,...$. We consider here only the case $m=8$.

Let us point out that the triangles of this lattice cannot be made
arbitrary small  in the standard Lobachevski metric (their size is
fixed in $H^2$ by the number $m$) but our domain $D$ can be made
arbitrarily large.

For example, for every finite set of vertices $K\subset H^2_m$ and
positive integer $r\in Z_+$ we define a domain $D_{K,r}$ consisting
of vertices $x\in D_{K,r}$ such that distance $d(x,K)$ is  no more
than $r$. We measure distances between vertices  (sets) counting
minimal number of edges needed for joining them. The simplest
important domains of that kind correspond to the cases: 1.$K$=vertex
$0$ (we call it ''standard ball'' $D_r$); 2.$K$=triangle $T$ i.e.
$K$=3 vertices of $T$. We denote it $D_{T,r}$. It is also like a
ball with center in the center of triangle. For $K=D_{K',l}$ we have
$D_{K,r}=D_{K',r+l}$ if $K$ is connected (i.e. no jumps of the
length more than one are needed to reach one point from another).

There is very big automorphism group mapping this Lattice into
itself but this group is non-commutative. It does not contain big
enough commutative subgroups, so nothing like Fourier transform
exists here. We cannot construct good enough Green function. Our
operators $Q^b,Q^w$ map space of functions of vertices into the
space of functions of black and white triangles correspondingly.
Therefore we cannot iterate them. So we cannot construct analogs of
polynomials here.

{\bf Problem: How to construct basis of d-holomorphic functions
$\psi_l(x)$ in $H^2_8$ which are globally bounded in all space?}

In the continuous case our space is realized as a unit disc
$D^2:|z|<1$. We have basis of bounded holomorphic functions $z^k,
k=0,1,2,...$. Every rational function with poles outside of unit
disc is bounded in it.

Easy to construct some d-holomorphic functions $z_{P,r}(x)$ equal to
zero inside of the $r$-ball $D_r$ (i.e. for $x\in D_l,l<r$) and
equal to zero along the boundary $\partial D_r$ except of the
specific place $P\subset \partial D_r$ looking like in Fig 6, a and
b:

\begin{center}
\color{black} Fig 6 \vspace{5mm}
\mbox{\epsfysize= 3cm\epsffile{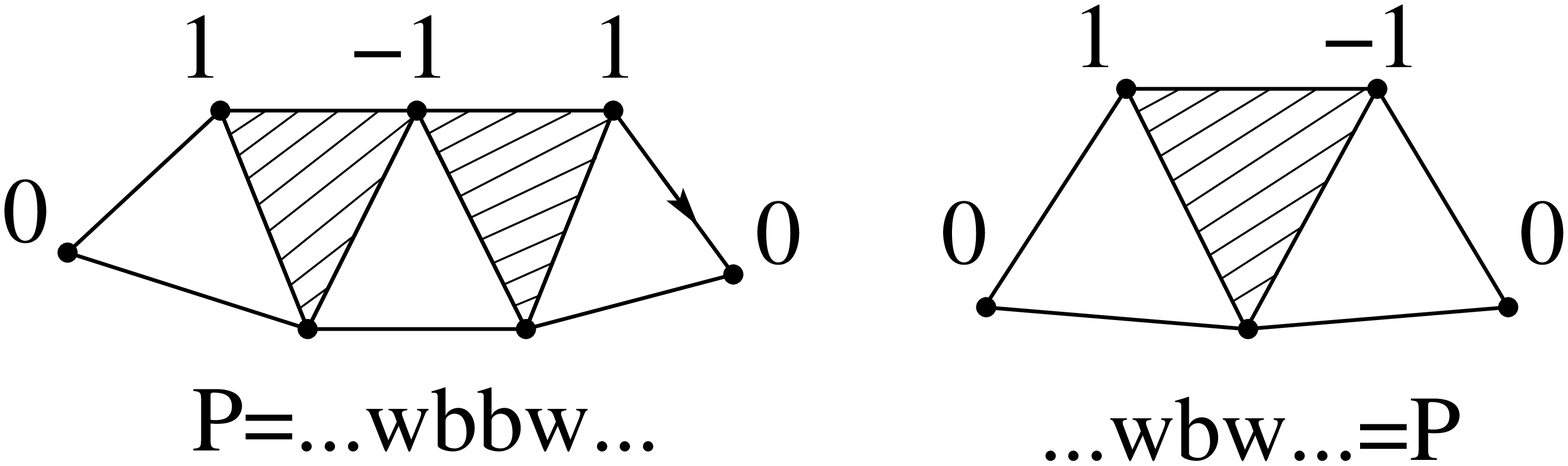}}
\end{center}

Its extension to the external domain is non-unique. {\bf Can we
construct a globally bounded extension?}.

Let us introduce class of right-convex paths.
\begin{defi} We call oriented path consisting of edges Right Convex
if it bounds two or three triangles only from the right site in
every vertex
\end{defi}

We are coding all right convex paths by the words in two symbols
$w,b$ assigning to every edge letter $w$ or $b$ depending on which
color has triangle from the right site of this path.-see Fig 7

\begin{center}\color{black}Fig 7\vspace{5mm}
\mbox{\epsfysize= 5cm\epsffile{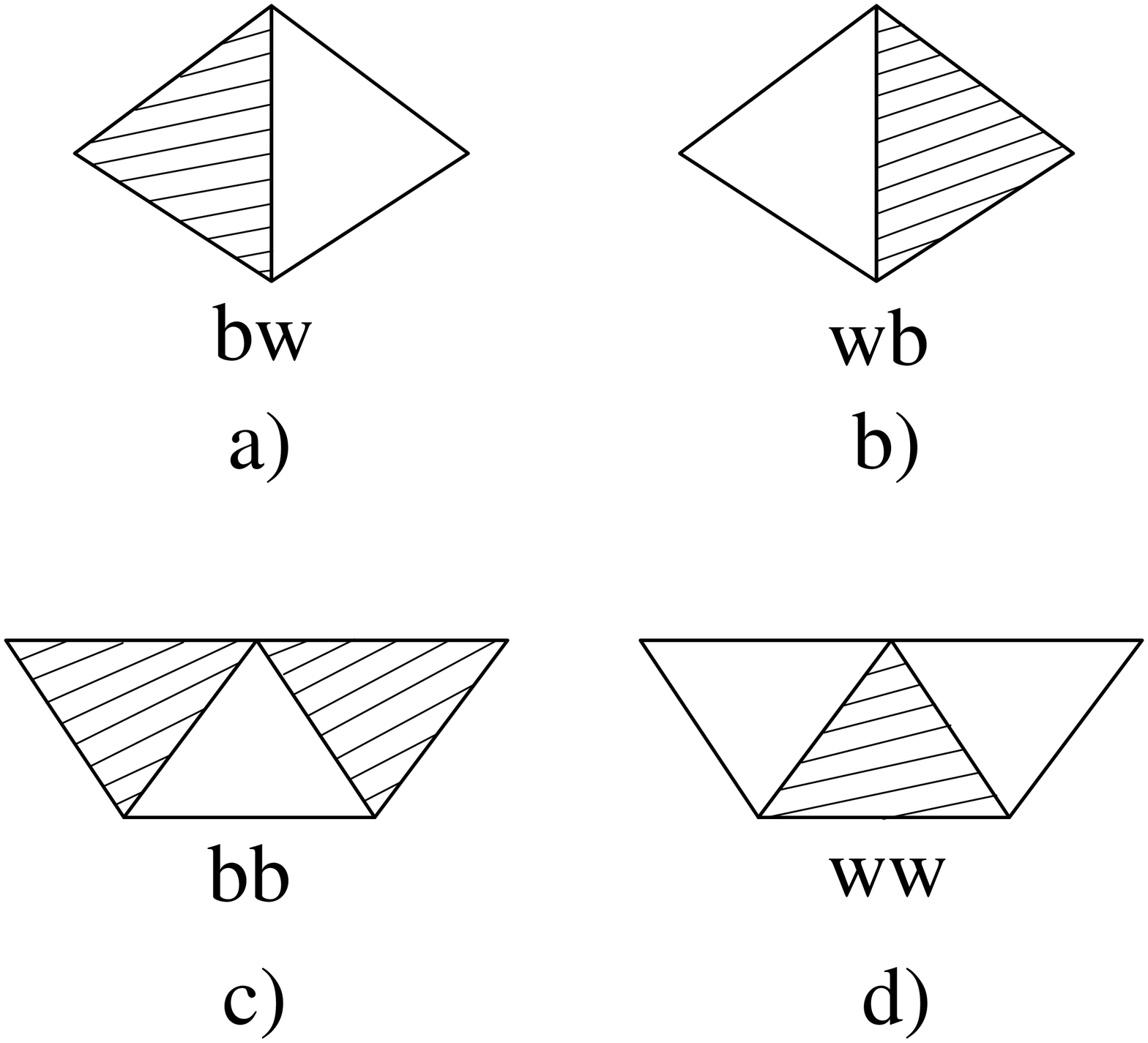}}
\end{center}

\begin{thm}
Let $\psi$ be a d-holomorphic function and $R$ is a full set of
zeroes $\psi(x)=0$. Consider  boundary of its complementary domain,
i.e. the set $D_{R,1} minus R$, which is a set of points-closest
neighbors of zeroes. Every connected component of this boundary set
is right convex choosing orientation such that corresponding
component of zeroes lies inside.
\end{thm}
Proof. In order to have zero value of $\psi$ in the point $x$ we
should have following values in the neighboring point outside of
zero set (see Fig 8, a and b).

\begin{center}\color{black}Fig 8\vspace{5mm}
\mbox{\epsfysize= 5cm\epsffile{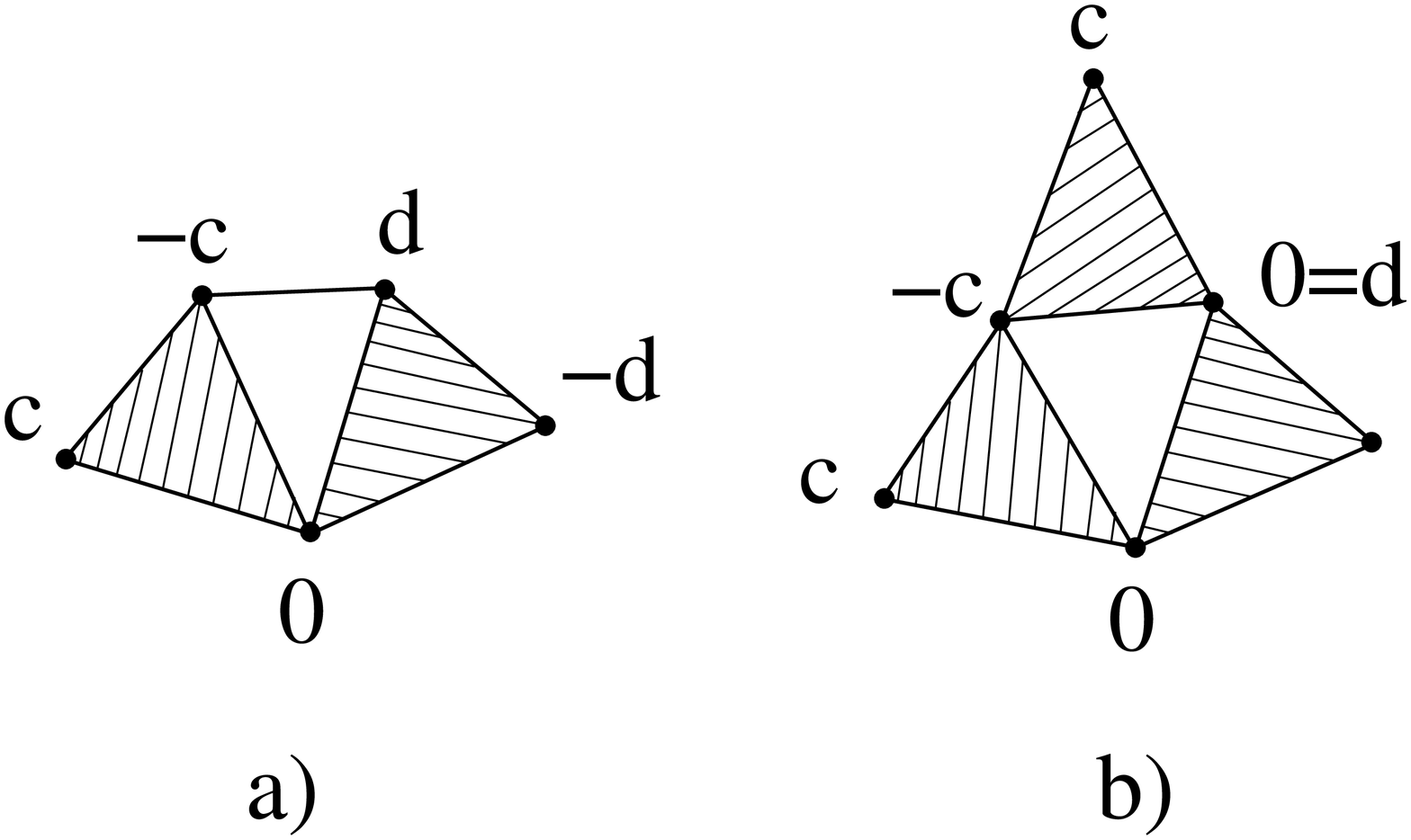}}
\end{center}

 No other possibilities exist because of
the equation $Q^b\psi=0$. It means in this case that $c\neq 0$
implies either $d\neq 0$ or $d=0$ in the Fig 8. In the first case
our boundary contains only two triangles inside. In the second case
$(d=0)$ our boundary contains three triangles inside. Anyway, it is
convex. Our theorem is proved.

\begin{center}
\color{black} Fig 9\vspace{5mm}
 \mbox{\epsfysize=5cm\epsffile{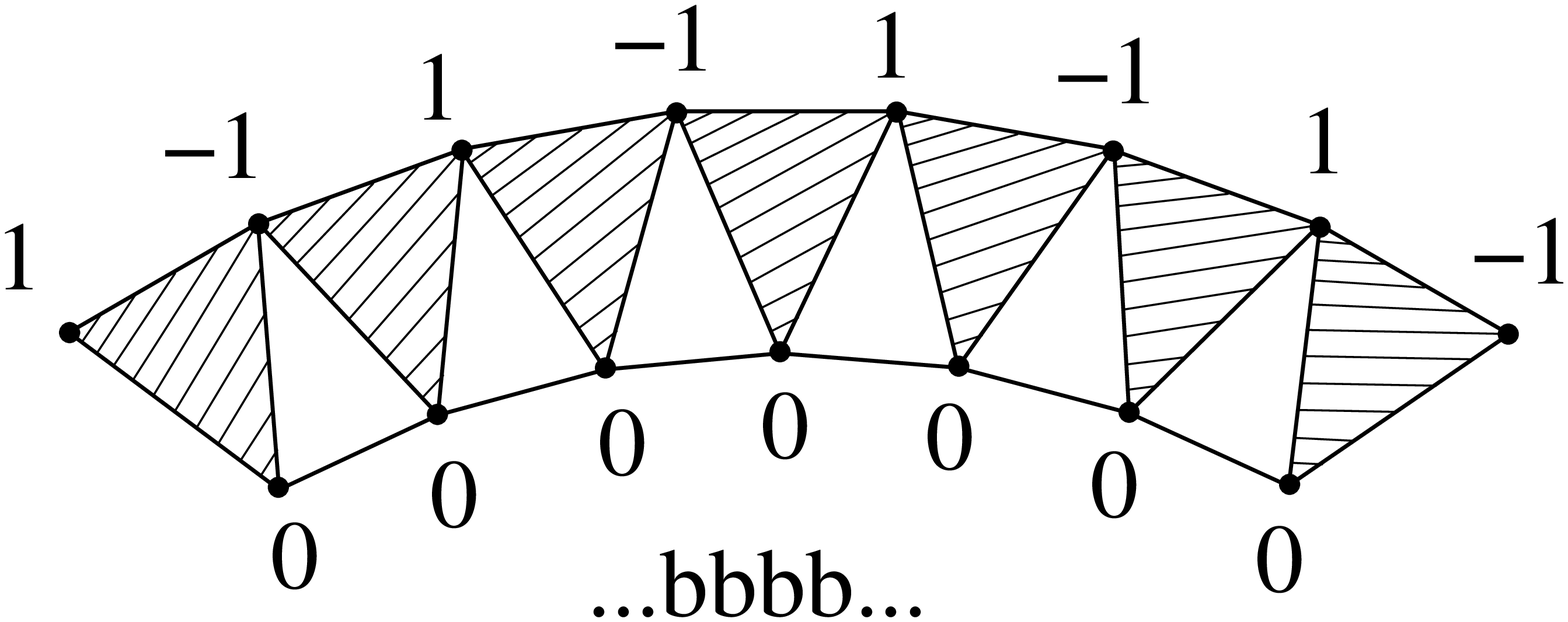}}
\end{center}

A maximal right convex set has a boundary path coded by the infinite
sequence $...bbbbb...$. For every point $x$ and direction-edge $l$
started in $x$ and having black triangle from the right, we uniquely
construct such maximal path $...bbb...$ and denote it by
$\gamma_{x,l}$. Easy to construct d-holomorphic function
$\psi_{x,l}(y)$ such that $\psi=0$ in the domain to the right of the
path $\gamma_{x,l}$, and $\psi(x)=\pm 1$ along the path
$\gamma_{x,l}$--see Fig 9. Its continuation to the complementary
domain is non-unique. We can define it from the requirement that the
growth is minimal if distance $d(y,\gamma_{x,l})\rightarrow \infty$,
but this definition is noneffective. {\bf How to find this growth?
What is a ''minimal'' function $\psi_{x,l}$?}. This function and its
group shifts give basis in the space of all d-holomorphic functions.

 Every
right-convex path $\gamma$ can be shifted to the left side by the
distance equal to one. We get a new  path $T(\gamma)$. Between these
two paths we have a strip containing all triangles $T$ having at
least one common vertex with path $\gamma$. New path $T(\gamma)$
consists of edges opposite to vertices belonging to $\gamma$. Only
triangles having exactly one common vertex with $\gamma$ participate
in this construction (see Fig 10).

\begin{center}
\color{black} Fig 10 \vspace{5mm}
\mbox{\epsfysize=5cm\epsffile{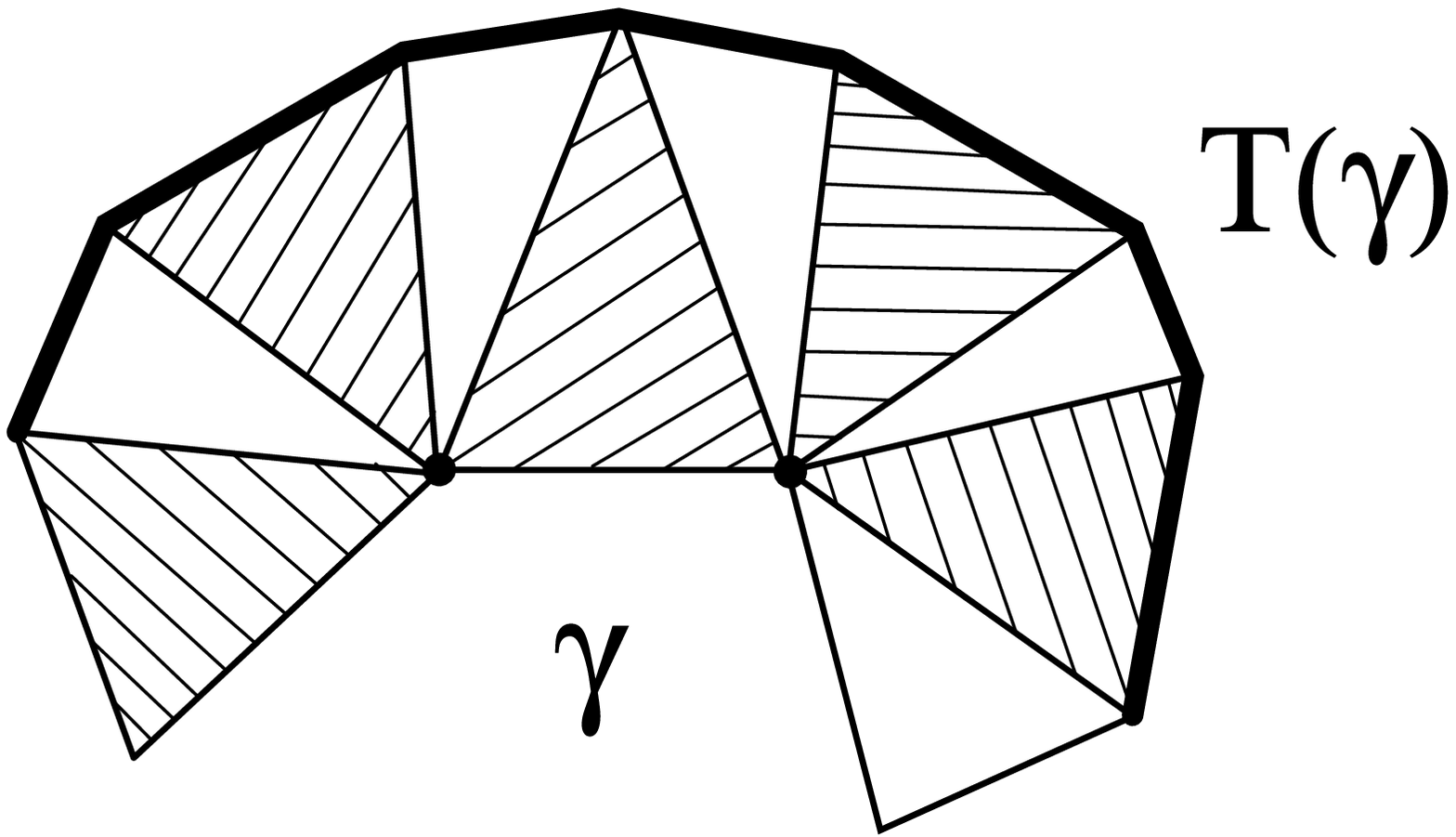}}
\end{center}

\begin{lem}
New path $T(\gamma)$ is also right convex. Its word  can be obtained
from the word describing $\gamma$ by the following procedure: take
every pair of neighboring letters in $\gamma$ and put between them
words written below $$bw\rightarrow bwbw, wb\rightarrow wbwb$$
$$bb\rightarrow bwb,ww\rightarrow wbw$$
After that delete old letters.What remains is exactly a word
$T(\gamma)$.
\end{lem}

Proof easily follows from the picture (see Fig 11).

\begin{center}\color{black} Fig 11
\vspace{5mm}
 \mbox{\epsfysize=10cm\epsffile{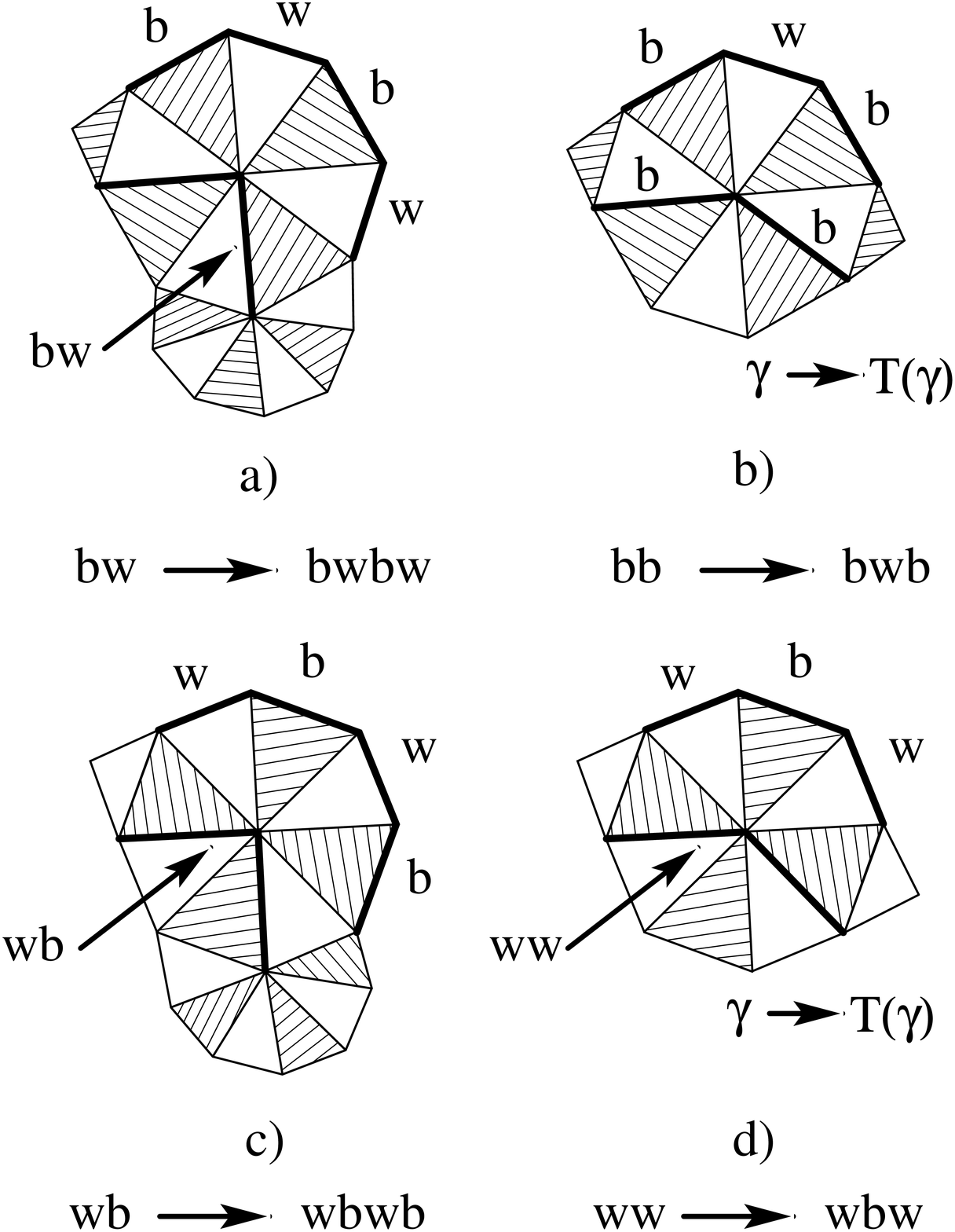}}
\end{center}

This map can be studied by the technic of symbolic dynamics. Mike
Boyle from the University of Maryland helped me a lot with this
business. We introduce new letters in order to describe this map.
Let $bw$ corresponds to the symbol $w_b$, $bb$ corresponds to $b_b$
$wb$ corresponds to $b_w$ and $ww$ corresponds to $w_w$. Written in
these symbols, our map $T$ has a form $$w_b\rightarrow w_bb_ww_bw_w,
w_w\rightarrow b_ww_bw_w$$ and $$b_w\rightarrow b_ww_bb_wb_b,
b_b\rightarrow w_bb_wb_b$$ After abelianization and replacing
product by sum, we are coming to the''Perron  matrix'' $A$ whose
largest eigenvalue is $\lambda=2+\sqrt{3}$. So we proved following

\begin{thm}The size  of  right convex path $\gamma$ shifted to the left side
by the distance one, increases asymptotically by the factor
$\lambda=2+\sqrt{3}$, so we have  $|T(\gamma)|\sim
(2+\sqrt{3})|\gamma|$. In particular, this is true for the boundary
of $r$-ball $|\gamma|=|\partial D_r|\sim \lambda^{r-k}|\partial
D_k|, k\geq 1$
\end{thm}

\begin{ex}We have $|\partial D_r|=8, 32, 120, 448, 1672, ...$ for
$r=1,2,3,4, 5,..$, so our asymptotic formula is practically exact
for $r\geq 3$.
\end{ex}

{\bf How to calculate dimension of the space of d-holomorphic
functions in the ball $D_r$? How many data on the boundary $\partial
D_r$ are needed to recover d-holomorphic function in $D_r$?}

Consider strip between $D_r$ and $D_{r+1}$-see Fig 12, r=1.

\begin{center}
\color{black} Fig 12 \vspace{5mm}
\mbox{\epsfysize=5cm\epsffile{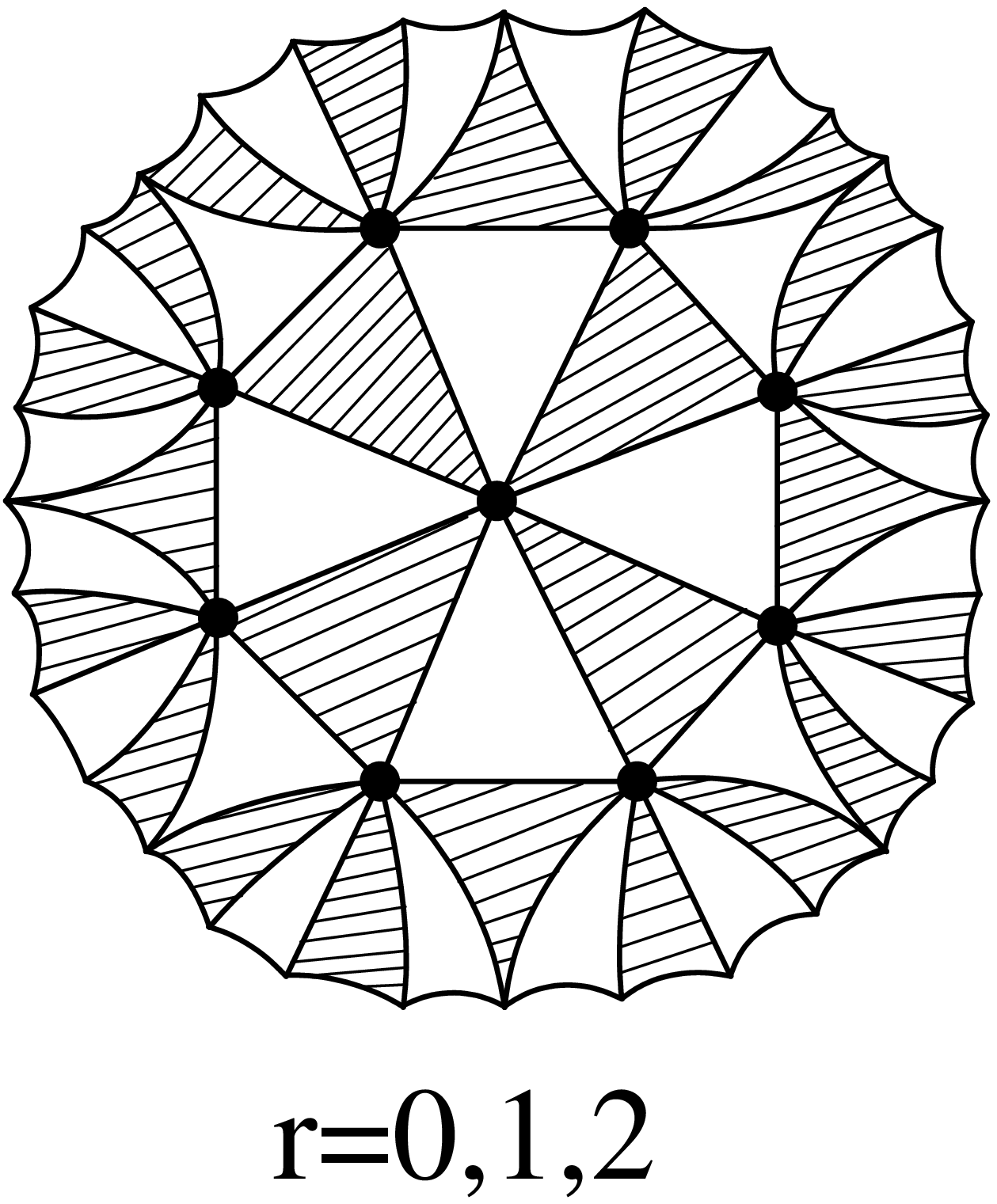}}
\end{center}

Every letter $b$ in $\partial D_{r+1}$ defines black triangle in the
strip touching $\partial D_r$ in one vertex. Every letter $w$ in
$\partial D_r$ defines black triangle in the same strip touching
$\partial D_{r+1}$ in one vertex. So total number of equations
$Q^b\psi=0$ in this strip is equal to $B_{r+1}+W_r$, the numbers of
black and white letters in the boundaries $\partial D_{r+1},\partial
D_r$. Unifying all strips with $k \leq r+1$, we get total number of
equations $Eq_{r+1}:Q^b\psi=0$:
$$Eq_{r+1}=(B_{r+1}+W_r)+...+(B_2+W_1)+B_1$$
For the number of points in $D_{r+1}$ we have
$$N_{r+1}=1+B_1+W_1+...+B_k+W_k+...+B_{r+1}+W_{r+1}$$
and $|\partial D_k|=B_k+W_k$. So we have $Eq_{r+1}=B_{r+1}+N_r-1$.
For the number of necessary data on the boundary $\partial D_{r+1}$
 we obtain (taking into account $W_k=B_k$)
 $$N_{r+1}-Eq_{r+1}=W_{r+1}+1=|\partial D_{r+1}|/2+1$$

{\bf Conclusion: Our result that ''The Number of Necessary Data is
equal to $[|\partial D_{r+1}|/2]+1$'' is an exact analog of the
Continuous case: Only the set of Fourier coefficients  on the
boundary circle corresponding to the exponents $\exp\{in\phi\}$ with
nonnegative $n\geq 0$, is needed for the reconstruction of
holomorphic function in the whole disc inside.}


\begin{thebibliography}{99}

\bibitem{DN}{I.Dynnikov, S.Novikov. Geometry of Triangle Equation, Moscow Math Journal-MMJ
(2003) v 3, pp 410-438 }


\bibitem{LF} {J.(Lelong)-Ferrand.  Fonctions preharmoniques et fonctions preholomorphes,
 Bull Sci Math(1944)  v 68 second series, pp 152-180}

\bibitem{DUF} {R.Duffin.  Basic properties of discrete analytic function, Duke Math Journal 1956 vol 23 pp 335-363}

\bibitem{N}{S.Novikov. Discrete $GL_n$-Connections. Proceeding of Steklov Math Institute,
(2004) v 247, pp 186-201}

\bibitem{Mer}{A.Bobenko, C.Mercat, Yu.Suris.  Linear and nonlinear theories of discrete
analytic functions. Integrable structure and isomonodromic Green functions,
 J.Reine Angew. Mathematics(2005) v 583, pp 117-161 }

\bibitem{DSS}{N.Dolbilin,  M.Shtan'ko, M.Shtogrin. a) Cubic subcomplexes in regular
lattices. Dokl. Akad. Nauk SSSR Dokl. 291, 277-279 (1986), English
translation: Sov. Math. Dokl. 34, 467-469 (1987); b) (jointly with
Sedrakyan, A.G.). A topology for the family of parametrizations of
two-dimensional cycles arising in the three-dimensional Ising model.
Dokl. Akad. Nauk SSSR  295, 12-23 (1987), English translation: Sov.
Math., Dokl. 36, No.1, 11-15 (1987); c) The problem of
parametrization of cycles modulo 2 in a three-dimensional
 cubic lattice. Izv. Akad. Nauk SSSR, Ser. Mat., 52, No. 2, 355-377 (1988);
 English translation:
 Math. USSR, Izv. 52, No.2, 359-383 (1989); ñ) Quadrillages and
parametrizations of lattice cycles. Tr. Mat. Inst. Steklova  196,
66-85 (1991); English translation: Proc. Steklov Inst. Math. 196,
73-93 (1992); d)  Cubic manifolds in lattices. Izv. Ross. Akad.
Nauk, Ser. Mat., 58, No. 2, 93-107 (1994), English translation:
Izv., Math. 44, No.2, 301-313 (1995). }

\bibitem{N2}{S.Novikov. Topology I. Encyclopedia Math Sciences (Editor General-R.Gamkrelidze), vol 13
 (S.Novikov, B.Botvinnik, R Burns--editors),
Springer Verlag (in english), page 41-42}

\bibitem{N1}{S.Novikov. Algebraic properties of 2D difference operators, Russian Math Surveys
(1997) v 52, n 1, pp 225-226 }

\bibitem{ND}{S.Novikov, I.Dynnikov. Discrete Spectral Symmetries of differential and difference
low dimensional operators, Russian Math Surveys (1997) v 52, n 5, pp
175-234}

\bibitem{GN}{ P.Grinevich, R.Novikov.  The Cauchy kernel for Novikov-Dynnikov (DN) discrete
complex analysis in triangular lattioces, Russian Math Surveys
(2007) v 62, n 4, pp 799-801}

\bibitem{GN1} {P.Grinevich, R.Novikov, to appear in 2010}

\end{thebibliography}
\end{document}